\documentclass[a4paper,12pt,reqno]{amsart}

\usepackage{amsmath}
\usepackage{amscd}
\usepackage{amssymb}
\usepackage{mathrsfs}
\usepackage[left=2.5cm,right=2.5cm,bottom=3cm,top=3cm]{geometry}
\newtheorem{thm}{Theorem}[section]

\newtheorem{prop}[thm]{Proposition}
\theoremstyle{definition}
\newtheorem{defn}[thm]{Definition}

\newtheorem{rem}[thm]{Remark}

\numberwithin{equation}{section}

\begin{document}

\title[On almost quasi-negative holomorphic sectional curvature]{On almost quasi-negative holomorphic sectional curvature}

\author{Yashan Zhang and Tao Zheng}

\thanks{Zhang is partially supported by Fundamental Research Funds for the Central Universities and National Natural Science Foundation of China (No. 12001179)}
\thanks{Zheng is partially supported by Beijing Institute of Technology Research Fund Program for Young Scholars}

\address{School of Mathematics, Hunan University, Changsha 410082, China}
\email{yashanzh@hnu.edu.cn}

\address{School of Mathematics and Statistics, Beijing Institute of Technology, 5 South Zhongguancun Street, Haidian, Beijing 100081, China}
\email{zhengtao08@amss.ac.cn}


\begin{abstract}
A recent celebrated theorem of Diverio-Trapani and Wu-Yau states that a compact K\"ahler manifold admitting a K\"ahler metric of quasi-negative holomorphic sectional curvature has an ample canonical line bundle, confirming a conjecture of Yau. In this paper we shall consider a natural notion of \emph{almost quasi-negative holomorphic sectional curvature} and extend this theorem to compact K\"ahler manifolds of almost quasi-negative holomorphic sectional curvature. We also obtain a gap-type theorem for the inequality $\int_Xc_1(K_X)^n>0$ in terms of the holomorphic sectional curvature. In the discussions, we introduce a capacity notion for the negative part of holomorphic sectional curvature, which plays a key role in studying the relation between the almost quasi-negative holomorphic sectional curvature and ampleness of the canonical line bundle. 
\end{abstract}

\maketitle

\section{Introduction}

\subsection{Backgrounds} Negativity of the holomorphic sectional curvature is a classic subject in the study of complex geometry, and deeply relates to several topics, including Kobayashi hyperbolicity, rigidity theorems of holomorphic maps, positivity of the canonical line bundle, etc.. One of the most fundamental questions is a conjecture of S.-T. Yau, which predicts that the negativity of the holomorphic sectional curvature should closely affect the positivity of the canonical line bundle. Precisely, Yau conjectured that \emph{a compact K\"ahler manifold admitting a K\"ahler metric of negative holomorphic sectional curvature has an ample canonical line bundle}, which has attracted lots of attention (see e.g. \cite{DT,HLW,HLW2,HLWZ,No,TY,W,WWY,WY1,WY2,YZ,Zys1} etc.). A recent breakthrough of Wu-Yau \cite{WY1} proved

\begin{thm}\cite{WY1}\label{thm1.1}
A projective manifold admitting a K\"ahler metric of negative holomorphic sectional curvature has an ample canonical line bundle.
\end{thm}
The surface case and threefold case in Theorem \ref{thm1.1} were previously proved by B. Wong \cite{W} and Heier-Lu-B. Wong \cite{HLW}, respectively, and  Heier-Lu-B. Wong \cite{HLW2} proved Theorem \ref{thm1.1} by assuming Abundance Conjecture.

Tosatti-Yang \cite{TY} extended Wu-Yau's Theorem \ref{thm1.1} to the K\"ahler case and hence proved above-mentioned Yau's conjecture in full generality:

\begin{thm}\cite{TY}\label{thm_ty'}
A compact K\"ahler manifold admitting a K\"ahler metric of negative holomorphic sectional curvature has an ample canonical line bundle.
\end{thm}

A key step in \cite{TY} is the following, whose projective case is a consequence of Yau's Schwarz Lemma \cite{Y78} and Mori's Cone Theorem (see e.g. \cite{Ma}):
\begin{thm}\cite{TY}\label{thm_ty}
A compact K\"ahler manifold admitting a K\"ahler metric of nonpositive holomorphic sectional curvature has a nef canonical line bundle.
\end{thm}
Wu-Yau and Tosatti-Yang's result was later extended to the case of quasi-negative holomorphic sectional curvature (which was also a conjecture of Yau) by Diverio-Trapani \cite{DT} and Wu-Yau \cite{WY2} (also see P. Wong-Wu-Yau \cite{WWY} for the special case that the manifold has Picard number one):
\begin{thm}\cite{DT,WY2}\label{thm_quasi}
A compact K\"ahler manifold admitting a K\"ahler metric of quasi-negative holomorphic sectional curvature has an ample canonical line bundle.
\end{thm}

In the proof of Theorem \ref{thm_quasi}, a key step is the following inequality (see \cite[Section 2]{DT} and \cite[Theorem 2]{WY2}):
\begin{thm}\cite{DT,WY2}\label{thm_quasi'}
Let $X$ be an $n$-dimensional compact K\"ahler manifold. If $X$ admits a K\"ahler metric of quasi-negative holomorphic sectional curvature, then
\begin{equation}
\int_Xc_1(K_X)^n>0\nonumber.
\end{equation}
\end{thm}

An excellent exposition on this topic can be found in a recent survey by Diverio \cite{Di}.
\begin{rem}
It is conjectured (see \cite[Remark 1.7]{TY}, \cite[Conjectures 1.1]{YZ}) that the above theorems could be extended to the Hermitian case, and there are interesting progresses, see \cite{YZ,Ta, LS}.
\end{rem}

\subsection{Almost nonpositivity of the holomorphic sectional curvature and nefness of the canonical line bundle} Motivated by Tosatti-Yang's Theorem \ref{thm_ty} and the definition of nefness of the canonical line bundle, the first-named author introduced in the previous work \cite{Zys1} a natural notion of \emph{almost nonpositive holomorphic sectional curvature}. Let's first recall some necessary notations. Let $(X,\omega)$ be a K\"ahler manifold of dimension $dim_{\mathbb C}X=n$. In a local holomorphic chart $(z^1,...,z^n)$, we write
$$\omega=\sqrt{-1}g_{i\bar j}dz^i\wedge d\bar z^j,$$
then the curvature tensor $R^\omega=\{R^\omega_{i\bar jk\bar l}\}$ of $\omega$ is given by
$$R^\omega_{i\bar jk\bar l}=-\frac{\partial^2g_{k\bar l}}{\partial z^i\partial\bar z^j}+g^{\bar qp}\frac{\partial g_{k\bar q}}{\partial z^i}\frac{\partial g_{p\bar l}}{\partial \bar z^j}.$$

Given $x\in X$ and $W\in T^{1,0}_xX\setminus\{0\}$, the \emph{holomorphic sectional curvature of $\omega$ at $x$ in the direction $W$} is
$$H^\omega_x(W):=\frac{R^\omega(W,\overline W,W,\overline W)}{|W|^4_\omega}.$$

We set
$$H^\omega_x:=\sup \{H^\omega_x(W)|W\in T^{1,0}_xX\setminus\{0\}\}$$
and
$$\mu_{\omega}:=\sup_{x\in X}H^\omega_x.$$
\emph{Throughout this paper, we always assume that each involved K\"ahler metric $\omega$ satisfies $\mu_{\omega}>0$.}\\

Define a continuous real function $\kappa_\omega$ for $(X,\omega)$ as follows:
\begin{equation}\label{kappa}
\kappa_\omega(x):=\rho(H^\omega_x)\cdot H^\omega_x,
\end{equation}
where $\rho:\mathbb R\to\{\frac{n+1}{2n},1\}$ is a function with $\rho(s)=\frac{n+1}{2n}$ for $s\le0$ and $\rho(s)=1$ for $s>0$.

\begin{defn}[Almost nonpositive holomorphic sectional curvature \cite{Zys1}]\label{defn1}
Let $(X,\omega_0)$ be a compact K\"ahler manifold.
\begin{itemize}
\item[(1)] Let $[\omega]$ be a K\"ahler class on $X$. Define a number  $\mu_{[\omega]}$ for $[\omega]$ in the following way: \\
\centerline{$\mu_{[\omega]}:=\inf\{\mu_{\omega'}|\omega'$ is a K\"ahler metric in $[\omega]\}$.}
\item[(2)] We say \textbf{$X$ is of almost nonpositive holomorphic sectional curvature} if there exist a sequence number $\epsilon_i\searrow0$ and a sequence of K\"ahler class $\alpha_i$ on $X$ such that $\mu_{\alpha_i}\alpha_i<\epsilon_i[\omega_0]$.
\item[(3)] We say \textbf{the K\"ahler class $[\omega]$ is of almost nonpositive holomorphic sectional curvature} if $\mu_{[\omega]}=0$.
\end{itemize}
\end{defn}

As we mentioned above, Definition \ref{defn1} is mainly motivated by the definition of the nefness of the canonical line bundle, which turns out to be also natural from the point of view of generalizing rigidity theorem for holomorphic maps, see \cite[Theorem 1.1, Corollary 4.1, etc.]{Zys2}. As an application, Theorem \ref{thm_ty} can be generalized as follows:
\begin{thm}\cite{Zys1}\label{thm_z1}
A compact K\"ahler manifold of almost nonpositive holomorphic sectional curvature has a nef canonical line bundle.
\end{thm}

For convenience let's extract the following particular case of Theorem \ref{thm_z1}:
\begin{thm}\cite{Zys1}\label{thm_z2}
A compact K\"ahler manifold admitting a K\"ahler class of almost nonpositive holomorphic sectional curvature has a nef canonical line bundle.
\end{thm}

\subsection{Motivation and main results}
Given Theorems \ref{thm_z1} and \ref{thm_z2}, it seems natural to ask: can we extend Theorem \ref{thm_quasi} to certain setting of ``\textbf{almost}" quasi-negative holomorphic sectional curvature? Here we propose a notion of almost quasi-negative holomorphic sectional curvature for a K\"ahler class, corresponding to Definition \ref{defn1} (3). Recall that on a compact K\"ahler manifold $X$, a K\"ahler class $[\omega_0]$ has almost nonpositive holomorphic sectional curvature, i.e. $\mu_{[\omega_0]}=0$, if and only if there exists a sequence of K\"ahler metrics $\hat\omega_i\in[\omega_0]$ on $X$ such that $\mu_{\hat\omega_i}\searrow0$. Then the following definition may be natural.

\begin{defn}[Almost quasi-negative holomorphic sectional curvature]\label{defn_quasi}
On a compact K\"ahler manifold $X$, we say \textbf{a K\"ahler class $[\omega_0]$ has almost quasi-negative holomorphic sectional curvature} if there exists a sequence of K\"ahler metrics $\hat\omega_i\in[\omega_0]$ on $X$ such that
\begin{itemize}
\item[(1)] $\mu_{\hat\omega_i}\searrow0$;
\item[(2)] $\kappa_{\hat\omega_i}<0$ somewhere on $X$.
\end{itemize}
\end{defn}

\begin{rem}
(1) If a K\"ahler metric $\omega_0$ has quasi-negative holomorphic sectional curvature, then its K\"ahler class $[\omega_0]$ obviously has almost quasi-negative holomorphic sectional curvature.\\
(2) A K\"ahler class of almost quasi-negative holomorphic sectional curvature must have almost nonpositive holomorphic sectional curvature.
\end{rem}

Given Theorems \ref{thm_quasi} and \ref{thm_z2}, it seems natural to expect that \emph{a compact K\"ahler manifold admitting a K\"ahler class of almost quasi-negative holomorphic sectional curvature has an ample canonical line bundle}. However, the condition in Definition \ref{defn_quasi} should be too weak to conclude the ampleness of the canonical line bundle, once we compare it with the Ricci curvature case, as explained in the following remark.

\begin{rem}\label{rem_ric}
We say a K\"ahler class $[\omega_0]$ has almost quasi-negative Ricci curvature if there is a number sequence $\epsilon_i\searrow0$ and a sequence of K\"ahler metrics $\hat\omega_i\in[\omega_0]$ such that
\begin{itemize}
\item[(1')] $Ric(\hat\omega_i)\le\epsilon_i\omega_0$;
\item[(2')] $Ric(\hat\omega_i)<0$ somewhere on $X$.
\end{itemize}
In general, such conditions (1') and (2') may not be enough to conclude 
\begin{equation}\label{c_1^n}
\int_X(2\pi c_1(K_X))^n=\int_X(-Ric(\hat\omega_i))^n>0,\,\,as\,\,i\to\infty,
\end{equation}
since the negative part of $Ric(\hat\omega_i)$ could be very small. To conclude \eqref{c_1^n}, it is natural to additionally require that, for example, there is a positive number $\delta$ such that for any $i$,
$$\int_{\{Ric(\hat\omega_i)<0\}}(-Ric(\hat\omega_i))^n\ge\delta,$$
which in some sense means that the ``capacity'' of the negative part of the Ricci curvature has a uniform positive lower bound and can help to control the small positive part of the Ricci curvature.
\end{rem}

Motivated by the above Remark \ref{rem_ric} on Ricci curvature case, to conclude the inequality $\int_X(c_1(K_X))^n>0$ for a compact K\"ahler manifold admitting a K\"ahler class of almost quasi-negative holomorphic sectional curvature, we may need to additionally require some uniform positive lower bounds on certain ``capacity'' of the negative part of the holomorphic sectional curvature. Note that the analog of $Ric(\hat\omega_i)$, which is invariant under rescaling the metric, for the holomorphic sectional curvature may be naturally chosen to be $\kappa_{\hat\omega_i}\hat\omega_i$. From this point of view, to conclude $\int_X(c_1(K_X))^n>0$ and the ampleness of the canonical line bundle under the conditions in Definition \ref{defn_quasi}, it is natural to additionally require that the negative part of $\kappa_{\hat\omega_i}\hat\omega_i$ has certain uniformly positive ``capacity'', which motivates us to consider the following definition.

\begin{defn}[Capacity of negative part of the holomorphic sectional curvature]\label{defn-cap} 
Let $X$ be an $n$-dimensional compact K\"ahler manifold, and $\omega_0$ a background K\"ahler metric on $X$. For any K\"ahler metric $\omega$ on $X$, $\lambda\in(0,\infty)$, we denote 
\begin{equation}\label{eq-U}
U(\omega,\lambda;\omega_0):=\{x\in X|\kappa_{\omega}<0\,\,and\,\,(-\kappa_{\omega}\omega)^n\le\lambda\omega_0^n\,\,at\,\,x\}
\end{equation}
and 
\begin{equation}\label{eq-V}
V(\omega,\lambda;\omega_0):=\{x\in X|\kappa_{\omega}<0\,\,and\,\,(-\kappa_{\omega}\omega)^n>\lambda\omega_0^n\,\,at\,\,x\},
\end{equation}
and define
\begin{equation}\label{eq-H}
\mathscr H(\omega,\lambda;\omega_0):=\int_{U(\omega,\lambda;\omega_0)}(-\kappa_{\omega}\omega)^n+\lambda\int_{V(\omega,\lambda;\omega_0)}\omega_0^n=\int_{\{\kappa_{\omega}<0\}}\min\left\{\frac{(-\kappa_{\omega}\omega)^n}{\lambda\omega_0^n},1\right\}\lambda\omega_0^n,
\end{equation}
we call $\mathscr H(\omega,\lambda;\omega_0)$ \textbf{the $\lambda$-capacity of negative part of the holomorphic sectional curvature of $\omega$ with respect to $\omega_0$}.
\end{defn}

Note the basic properties of $\mathscr H(\omega,\lambda;\omega_0)$:
\begin{itemize}
\item[(1)] For any K\"ahler metric $\omega$ with $\{\kappa_\omega<0\}\neq\emptyset$, we have $\mathscr H(\omega,\lambda;\omega_0)>0$ for each $\lambda\in(0,\infty)$;
\item[(2)] $\mathscr H(\omega,\lambda;\omega_0)$ is incresing in $\lambda\in(0,\infty)$;
\item[(3)] For any fixed $\omega$ and $\omega_0$, there is a positive number $\Lambda$, which depends on $\omega$ and $\omega_0$, such that $\frac{(-\kappa_{\omega}\omega)^n}{\Lambda\omega_0^n}\le1$ on $\{\kappa_{\omega}<0\}$, and hence
$\mathscr H(\omega,\lambda;\omega_0)\equiv\mathscr H(\omega,\Lambda;\omega_0)=\int_{\{\kappa_{\omega}<0\}}(-\kappa_{\omega}\omega)^n$ for any $\lambda\ge\Lambda$. Therefore, $\mathscr H(\omega,\lambda;\omega_0)$ is incresing in $\lambda\in(0,\infty)$, and always eventurally stabilizes to $\int_{\{\kappa_{\omega}<0\}}(-\kappa_{\omega}\omega)^n$ as $\lambda\to\infty$.
\item[(4)] In the above definition, the reference metric $\omega_0$ is involved. However, in the followings we will focus on the case that the capacities for a sequence of metrics admit a uniform positive lower bound, which does not depend on the choice of the reference metric, see discussions after the statement of Theorem \ref{thm_main}.
\end{itemize}

Comparing with the Ricci curvature case described in Remark  \ref{rem_ric}, it is natural to expect that the above capacity notions shall play a role in studying the relation between the almost quasi-negative holomorphic sectional curvature and ampleness of the canonical line bundle. 

The following is the main result of this paper, generalizing Diverio-Trapani and Wu-Yau's Theorem \ref{thm_quasi} and indicating the importance of the above capacity notion.
\begin{thm}\label{thm_main}
Let $X$ be an $n$-dimensional compact K\"ahler manifold, $\omega_0$ a K\"ahler metric on $X$ and assume that $[\omega_0]$ has almost quasi-negative holomorphic sectional curvature. Then we fix a  family of K\"ahler metrics $\hat\omega_i\in[\omega_0]$ satisfying conditions (1) and (2) in Definition \ref{defn_quasi}. Assume additionally that there is some $\lambda_0\in(0,\infty)$ such that
\begin{equation}\label{ass'n0} 
\limsup_{i\to\infty}\mathscr H(\hat\omega_i,\lambda_0;\omega_0)>0.
\end{equation}
Then the canonical line bundle $K_X$ of $X$ is ample.
\end{thm}

In Theorem \ref{thm_main}, up to passing to a subsequence, assumption \eqref{ass'n0} is equivalent to that one of the followings holds:
\begin{itemize}
\item[(i)]
\begin{equation}\label{ass'n}
\limsup_{i\to\infty}\int_{U(\hat\omega_i,\lambda_0;\omega_0)}(-\kappa_{\hat\omega_i}\hat\omega_i)^n>0,
\end{equation}
\item[(ii)]
\begin{equation}\label{ass'n'}
\limsup_{i\to\infty}\int_{V(\hat\omega_i,\lambda_0;\omega_0)}\omega_0^n>0.
\end{equation}
\end{itemize}
Using this equivalence, it is easy to see that the assumption \eqref{ass'n0}  does not depend on the choice of the reference metric $\omega_0$. More precisely, if \eqref{ass'n0} holds for some $\omega_0$ and $\lambda_0$, then for any other K\"ahler metric $\omega_1$ on $X$, we can find a positive number $\lambda_1$ such that $\limsup_{i\to\infty}\mathscr H(\hat\omega_i,\lambda_1;\omega_1)>0.$

Moreover, the \eqref{ass'n'} somehow means that the points at which the holomorphic sectional curvature of $\hat\omega_i$ is ``heavily negative'' (i.e. $(-\kappa_{\hat\omega_i}\hat\omega_i)^n>\lambda\omega_0^n$) do not consentrate on a zero-measure set, as $i\to\infty$. Indeed, when the heavily negative part of the holomorphic sectional curvature has certainly large capacity, the \eqref{ass'n'} is automatically true, as our next result indicates.

\begin{thm}\label{thm_main'}
Let $X$ be an $n$-dimensional compact K\"ahler manifold, $\omega_0$ a K\"ahler metric on $X$ and assume that $[\omega_0]$ has almost quasi-negative holomorphic sectional curvature. Then we fix a  family of K\"ahler metrics $\hat\omega_i\in[\omega_0]$ satisfying conditions (1) and (2) in Definition \ref{defn_quasi}. There exists a number $K_0\ge1$, depending only on $\omega_0$, such that if $\lambda\in[1,\infty)$ satisfies
\begin{equation}\label{ass'n''} 
\lim_{i\to\infty}\int_{V(\hat\omega_i,\lambda;\omega_0)}\left(\log\frac{(-\kappa_{\hat\omega_i}\hat\omega_i)^n}{\omega_0^n}\right)\omega_0^n>K_0,
\end{equation}
then
$$\liminf_{i\to\infty}\int_{V(\hat\omega_i,\lambda;\omega_0)}\omega_0^n>0,$$ 
and hence the canonical line bundle $K_X$ of $X$ is ample.
\end{thm}

Theorem \ref{thm_main'} in particular shows that a certainly large capacity of the negative part of the holomorphic sectional curvature actually can NOT consentrate on a set of very small measure.

\begin{rem}
\begin{itemize}
\item[(1)] We may mention that in the above theorems we do not impose any conditions on the involved $\hat\omega_{i}$'s on the level of \emph{metrics}. 
\item[(2)]  Consequently, a compact K\"ahler manifold satisfying the assumptions in Theorems \ref{thm_main} or \ref{thm_main'} admits a K\"ahler-Einstein metric of negative scalar curvature, thanks to the fundamental theorem of Aubin \cite{A} and Yau \cite{Y}.
\end{itemize}
\end{rem}

\subsubsection{A gap theorem for $\int_Xc_1(K_X)^n>0$}
As in Diverio-Trapani and Wu-Yau's works \cite{DT,WY2}, a key step to prove our main Theorems \ref{thm_main} and \ref{thm_main'} is to obtain the inequality $\int_Xc_1(K_X)^n>0$. Note that the inequality $\int_Xc_1(K_X)^n>0$ is an \emph{openness condition}. In particular, given Theorem \ref{thm_quasi'}, it seems natural to expect that if a compact K\"ahler manifold admits a K\"ahler metric whose holomorphic sectional curvature is ``sufficiently approximate" to being quasi-negative, then $\int_Xc_1(K_X)^n>0$ should hold. From this point of view, it is very natural to explore gap-type theorems for $\int_Xc_1(K_X)^n>0$ in terms of holomorphic sectional curvature. Regarding this, we shall prove the following, which generalizes the above Theorem \ref{thm_quasi'} and may be seen as a gap theorem for $\int_Xc_1(K_X)^n>0$. 

\begin{thm}\label{thm_gap}
Let $X$ be an $n$-dimensional compact K\"ahler manifold. Arbitrarily fix a K\"ahler metric $\omega_0$ on $X$ and two positive constants $\delta_1,\delta_2$. Then there exists a sufficiently small positive number $\hat \epsilon$, which depends on the given $\omega_0$ and $\delta_1,\delta_2$, satisfying the followings. Assume that there exists a K\"ahler metric $\hat\omega\in[\omega_0]$ such that
\begin{itemize}
\item[(a)] $\mu_{\hat\omega}\le\hat\epsilon$ on $X$;
\item[(b)] $\mathscr H(\hat\omega,\delta_1;\omega_0)\ge\delta_2.$
\end{itemize}
Then
\begin{equation}\label{key_ineq}
\int_Xc_1(K_X)^n>0.
\end{equation}
\end{thm}

\begin{rem}\label{rem-gap}
Under the assumptions in Theorem \ref{thm_gap}, actually the integral $\int_Xc_1(K_X)^n$ has a positive lower bound depending only on $\omega_0, \delta_1,\delta_2$, see Section \ref{pf_main}.
\end{rem}

\subsection{Organization} After recalling some necessary facts in the next section, we will prove Theorems \ref{thm_main}, \ref{thm_main'} and \ref{thm_gap} in Section \ref{pf_main}.

\section{Preliminaries}
Let's firstly collect several fundamental results in K\"ahler geometry, which will be used in the proofs of our theorems.

\subsection{Schwarz Lemma}
We will need the following Schwarz Lemma, see \cite{Y,Y78,R,WWY,WY1} (also \cite{Zys1}).
\begin{prop}[Schwarz Lemma]\label{pre1}
Let $\omega,\hat\omega$ be two K\"ahler metrics on $X$. Assume 
$$Ric(\omega)\ge-\omega+s\hat\omega,$$
where $s$ is a positive number. Then
$$\Delta_{\omega}\log tr_{\omega}\hat\omega\ge\left(-\kappa_{\hat\omega}+\frac{s}{n}\right)tr_{\omega}\hat\omega-1.$$
\end{prop}
Here, $\kappa_{\hat\omega}$ is the continuous function defined in \eqref{kappa}.

\subsection{Tian's $\alpha$-invariant}
Tian's $\alpha$-invariant \cite{Ti} is a key invariant in several complex/algebraic geometry problems, which will also play a particularly important role in our discussions.
\begin{prop}[Tian's $\alpha$-invariant]\label{pre2}
Let $\omega_0$ be a K\"ahler metric on $X$. Then there is a positive number $\alpha:=\alpha(X,[\omega_0])$ which depends only on $X$ and the K\"ahler class $[\omega_0]$ and satisfies the followings: for any $\beta\in(0,\alpha)$, we can find a positive number $C_{\beta}$ depending only on $\beta,\omega_0$ such that for any $u\in C^\infty(X)$ with $\omega_0+\sqrt{-1}\partial\bar\partial u>0$ and $\sup_Xu=0$,
$$\int_Xe^{-\beta u}\omega_0^n\le C_\beta.$$
\end{prop}

Note that for a positive number $c$, $\alpha(X,[c\omega_0])=c^{-1}\alpha(X,[\omega_0])$, and hence when $c$ is sufficiently small, the $\alpha$-invariant of $[c\omega_0]$ could be sufficiently large.

\subsection{Hartogs Lemma}
We will also need the following Hartogs-type lemma proved in \cite[Proposition 2.7]{GZ}.
\begin{prop}\label{pre3}
Let $\omega_0$ be a K\"ahler metric on $X$. Then there is a positive number $C$ which depends only on $X$ and $\omega_0$ and satisfies the followings: for any $u\in C^\infty(X)$ with $\omega_0+\sqrt{-1}\partial\bar\partial u>0$ and $\sup_Xu=0$,
$$\int_Xu\omega_0^n\ge -C.$$
\end{prop}

\section{Proofs of theorems}\label{pf_main}
We will firstly prove the gap Theorem \ref{thm_gap}, and then Theorem \ref{thm_main}, and finally Theorem \ref{thm_main'}.

\subsection{Proof of Theorem \ref{thm_gap}}\label{proof_gapthm} The arguments will make use of the similar strategy of Diverio-Trapani \cite{DT}. However, in our case, as the given K\"ahler metrics  no longer have nonpositively signed curvature in the pointwise sense, we have to make more efforts to control the extra terms caused by the positive part of the curvature, in which the key roles are played by Tian's $\alpha$-invariant and our capacity notion in Definition \ref{defn-cap}.

\subsubsection{Wu-Yau's continuity method}
Consider Wu-Yau's continuity method of $\varphi(t)$:
\begin{equation}\label{cm}
(t\hat\omega-Ric(\omega_0)+\sqrt{-1}\partial\bar\partial\varphi(t))^n=e^{\varphi(t)}\omega_0^n,
\end{equation}
which, writing $\omega(t):=t\hat\omega-Ric(\omega_0)+\sqrt{-1}\partial\bar\partial\varphi(t)$, is equivalent to
\begin{equation}\label{cm1}
Ric(\omega(t))=-\omega(t)+t\hat\omega.
\end{equation}
By the Schwarz Lemma (see Propopsition \ref{pre1}) we know
\begin{equation}\label{schlem}
\Delta_{\omega(t)}\log tr_{\omega(t)}\hat\omega\ge(-\kappa_{\hat\omega}+\frac{t}{n})tr_{\omega(t)}\hat\omega-1.
\end{equation}
Consequently the continuity method can be smoothly solved for $t\in(n\mu_{\hat\omega},+\infty)$, see \cite{TY} (or \cite{Zys1}).

Firstly, we fix a nonnegative number $b_0$ such that
\begin{equation}\label{riclower}
Ric(\omega_0)\ge-b_0\omega_0\,\,on\,\,X,
\end{equation}
which only depends on $\omega_0$.
Set
$$\Phi(t):=\varphi(t)+t\psi,$$
where $\psi\in C^\infty(X)$ satisfies $\hat\omega=\omega_0+\sqrt{-1}\partial\bar\partial\psi$ and $\sup_X\psi=0$.
Then the continuity method \eqref{cm} is equivalent to
\begin{equation}\label{cm-1}
(t\omega_0-Ric(\omega_0)+\sqrt{-1}\partial\bar\partial\Phi(t))^n=e^{\Phi(t)-t\psi}\omega_0^n,
\end{equation}

Using the maximum principle and the fact that $\psi\le0$ gives
\begin{equation}\label{phiupper}
\sup_X\Phi(t)\le\log(2n\hat\epsilon+b_0)^n\le\log(c_0+b_0)^n
\end{equation}
for $t\in(n\mu_{\hat\omega},2n\mu_{\hat\omega}]$, where $c_0$ is a positive number such that the $\alpha$-invariant $\alpha(X,c_0\omega_0)\ge2$, and we choose $\hat\epsilon\le \frac{c_0}{2n}$.

\subsubsection{A positive lower bound for quotients $\frac{-\int_{\{\kappa_{\hat\omega}<0\}}\kappa_{\hat\omega}e^{\frac{1}{n}t\psi}\left(\frac{\hat\omega^n}{\omega_0^n}\right)^{\frac1n}\omega(t)^n}{\int_X\omega(t)^n}$}

Next we  look at the quotient
\begin{align}\label{quotient'}
\frac{-\int_{\{\kappa_{\hat\omega}<0\}}\kappa_{\hat\omega}e^{\frac{1}{n}t\psi}\left(\frac{\hat\omega^n}{\omega_0^n}\right)^{\frac1n}\omega(t)^n}{\int_X\omega(t)^n} &=\frac{-\int_{\{\kappa_{\hat\omega}<0\}}\kappa_{\hat\omega}e^{\Phi(t)}e^{-(1-\frac1n)t\psi}\left(\frac{\hat\omega^n}{\omega_0^n}\right)^{\frac1n}\omega_0^n}
{\int_Xe^{\Phi(t)}e^{-t\psi}\omega_0^n} \nonumber\\
&=\frac{-\int_{\{\kappa_{\hat\omega}<0\}}\kappa_{\hat\omega}e^{\Phi(t)^*}e^{-(1-\frac1n)t\psi}\left(\frac{\hat\omega^n}{\omega_0^n}\right)^{\frac1n}\omega_0^n}{\int_Xe^{\Phi(t)^*}e^{-t\psi}\omega_0^n},
\end{align}
where $\Phi(t)^*:=\Phi(t)-\sup_X\Phi(t)$. We want to estimate this quotient from below by a positive number.

Firstly, we estimate $\int_Xe^{\Phi(t)^*}e^{-t\psi}\omega_0^n$ from above using $\alpha$-invariant. Indeed, since $0<t\le2n\mu_{\hat\omega}\le c_0$, we have
\begin{align}
0< t(\omega_0+\sqrt{-1}\partial\bar\partial\psi)\le c_0\omega_0+\sqrt{-1}\partial\bar\partial(t\psi)\nonumber
\end{align}
i.e. $t\psi\in PSH(X,c_0\omega_0)$ for $t\in(n\mu_{\hat\omega},2n\mu_{\hat\omega}]$. By the choice of $c_0$ (i.e. $\alpha(X,c_0\omega_0)\ge2$), and Propopsition \ref{pre2},
\begin{align}\label{eq_upper}
\int_Xe^{\Phi(t)^*}e^{-t\psi}\omega_0^n&\le\int_Xe^{-t\psi}\omega_0^n 
 =c_0^{-n}\int_Xe^{-t\psi}(c_0\omega_0)^n 
 \le c_0^{-n}C_{c_0\omega_0},
\end{align}
where $C_{c_0\omega_0}$ is a positive number depending only on $c_0\omega_0$. 

Secondly, we estimate $-\int_{\{\kappa_{\hat\omega}<0\}}\kappa_{\hat\omega}e^{\Phi(t)^*}e^{-(1-\frac1n)t\psi}\left(\frac{\hat\omega^n}{\omega_0^n}\right)^{\frac1n}\omega_0^n$ from below using the capacity assumption. Recall that for $t\in(n\mu_{\hat\omega},2n\mu_{\hat\omega}]$, $\Phi^*(t)\in PSH(X,(c_0+b_0)\omega_0)$ and $\sup_X\Phi^*(t)=0$. Then by Propopsition \ref{pre3} there exists a positive number $c_1$ depending only on $(c_0+b_0)\omega_0$ and hence only on $\omega_0$ such that
\begin{equation}\label{intlower}
\int_X\Phi^*(t)\omega_0^n\ge-c_1.
\end{equation}

Also, we obviously have $\{\kappa_{\hat\omega}<0\}=U(\hat\omega,\delta_1;\omega_0)\cup V(\hat\omega,\delta_1;\omega_0)=:\hat U\cup \hat V$, where $\hat U:=U(\hat\omega,\delta_1;\omega_0)$ and $\hat V:=V(\hat\omega,\delta_1;\omega_0)$, and so
\begin{align}\label{est_lower}
&-\int_{\{\kappa_{\hat\omega}<0\}}\kappa_{\hat\omega}e^{\Phi(t)^*}e^{-(1-\frac1n)t\psi}\left(\frac{\hat\omega^n}{\omega_0^n}\right)^{\frac1n}\omega_0^n\nonumber\\
&=-\int_{\hat U}\kappa_{\hat\omega}e^{\Phi(t)^*}e^{-(1-\frac1n)t\psi}\left(\frac{\hat\omega^n}{\omega_0^n}\right)^{\frac1n}\omega_0^n-\int_{\hat V}\kappa_{\hat\omega}e^{\Phi(t)^*}e^{-(1-\frac1n)t\psi}\left(\frac{\hat\omega^n}{\omega_0^n}\right)^{\frac1n}\omega_0^n\nonumber\\
&=:I+II.
\end{align} 
Applying Jenson inequality, \eqref{intlower} and the facts that $\Phi^*(t)\le0$ and $(-\kappa_{\hat\omega}\hat\omega)^n\le\delta_1\omega_0^n$ on $\hat U$, we have
\begin{align}\label{est_I}
I:&=-\int_{\hat U}\kappa_{\hat\omega}e^{\Phi(t)^*}e^{-(1-\frac1n)t\psi}\left(\frac{\hat\omega^n}{\omega_0^n}\right)^{\frac1n}\omega_0^n\nonumber\\
&=\int_{\hat U}e^{\Phi(t)^*}e^{-(1-\frac1n)t\psi}\left(\frac{(-\kappa_{\hat\omega}\hat\omega)^n}{\omega_0^n}\right)^{\frac1n-1}(-\kappa_{\hat\omega}\hat\omega)^n\nonumber\\
&\ge\delta_1^{\frac1n-1}\int_{\hat U}e^{\Phi(t)^*}(-\kappa_{\hat\omega}\hat\omega)^n\nonumber\\
&\ge\delta_1^{\frac1n-1}\left(\int_{\hat U}(-\kappa_{\hat\omega}\hat\omega)^n\right)\exp\left\{\frac{1}{\int_{\hat U}(-\kappa_{\hat\omega}\hat\omega)^n}\int_{\hat U}\Phi^*(t)(-\kappa_{\hat\omega}\hat\omega)^n\right\}\nonumber\\
&\ge\delta_1^{\frac1n-1}\left(\int_{\hat U}(-\kappa_{\hat\omega}\hat\omega)^n\right)\exp\left\{\frac{1}{\int_{\hat U}(-\kappa_{\hat\omega}\hat\omega)^n}\int_{\hat U}\Phi^*(t)(\delta_1\omega_0^n)\right\}\nonumber\\
&\ge\delta_1^{\frac1n-1}\left(\int_{\hat U}(-\kappa_{\hat\omega}\hat\omega)^n\right)\exp\left\{\frac{1}{\int_{\hat U}(-\kappa_{\hat\omega}\hat\omega)^n}\int_{X}\Phi^*(t)(\delta_1\omega_0^n)\right\}\nonumber\\
&\ge\delta_1^{\frac1n-1}\left(\int_{\hat U}(-\kappa_{\hat\omega}\hat\omega)^n\right)\exp\left\{\frac{-\delta_1c_1}{\int_{\hat U}(-\kappa_{\hat\omega}\hat\omega)^n}\right\}.
\end{align}

On the other hand, we can also estimate the second term as follows:
\begin{align}\label{est_II}
&II:=-\int_{\hat V}\kappa_{\hat\omega}e^{\Phi(t)^*}e^{-(1-\frac1n)t\psi}\left(\frac{\hat\omega^n}{\omega_0^n}\right)^{\frac1n}\omega_0^n\nonumber\\
&\ge\left(\int_{\hat V}\omega_0^n\right)\exp\left(\frac{1}{\int_{\hat V}\omega_0^n}\int_{\hat V}\log\left(-\kappa_{\hat\omega}\left(\frac{\hat\omega^n}{\omega_0^n}\right)^{1/n}\right)\omega_0^n+\frac{1}{\int_{\hat V}\omega_0^n}\int_{\hat V}\Phi^*(t)\omega_0^n\right)\nonumber\\
&\ge\left(\int_{\hat V}\omega_0^n\right)\exp\left(\frac{1}{n\int_{\hat V}\omega_0^n}\int_{\hat V}\left(\log\frac{(-\kappa_{\hat\omega}\hat\omega)^n}{\omega_0^n}\right)\omega_0^n-\frac{c_1}{\int_{\hat V}\omega_0^n}\right)\nonumber\\
&\ge\left(\int_{\hat V}\omega_0^n\right)\exp\left(\frac{\log\delta_1}{n}-\frac{c_1}{\int_{\hat V}\omega_0^n}\right)\nonumber\\
&=\delta_1^{1/n}\left(\int_{\hat V}\omega_0^n\right)\exp\left(-\frac{c_1}{\int_{\hat V}\omega_0^n}\right).
\end{align}

Now we apply conditions (b), i.e. 
$$\int_{\hat U}(-\kappa_{\hat\omega}\hat\omega)^n+\delta_1\int_{\hat V}\omega_0^n\ge\delta_2,$$ 
to see that one of the followings holds:
\begin{itemize}
\item[(b.1)] $\int_{\hat U}(-\kappa_{\hat\omega}\hat\omega)^n\ge\delta_2/2$; \item[(b.2)] $\delta_1\int_{\hat V}\omega_0^n\ge\delta_2/2$. 
\end{itemize}
In the case (b.1), combining $\int_{\hat U}(-\kappa_{\hat\omega}\hat\omega)^n\le\delta_1\int_{\hat U}\omega_0^n\le\delta_1\int_X\omega_0^n$, we use \eqref{est_I} to see that
\begin{align}\label{est_I'}
I\ge c_2:=\inf\left\{\delta_1^{\frac1n-1}s\exp(-\delta_1c_1/s)|\delta_2/2\le s\le\delta_1\int_X\omega_0^n\right\},
\end{align}
in which $c_2$ obviously depends only on $\omega_0,\delta_1,\delta_2$; while in the case (b.2), we use \eqref{est_II} to see that
\begin{align}\label{est_II'}
II\ge c_2':=\inf\left\{\delta_1^{1/n}s\exp(-c_1/s)|\frac{\delta_2}{2\delta_1}\le s\le\int_X\omega_0^n\right\},
\end{align}
in which $c_2'$ depends only on $\omega_0,\delta_1,\delta_2$.

In conclusion, combining 
\eqref{est_lower}, \eqref{est_I}, \eqref{est_II}, \eqref{est_I'} and \eqref{est_II'} gives
\begin{align}\label{est_lower'}
&-\int_{\{\kappa_{\hat\omega}<0\}}\kappa_{\hat\omega}e^{\Phi(t)^*}e^{-(1-\frac1n)t\psi}\left(\frac{\hat\omega^n}{\omega_0^n}\right)^{\frac1n}\omega_0^n\ge\min\{c_2,c_2'\}
\end{align}

Given \eqref{eq_upper}, \eqref{est_lower'} and \eqref{quotient'}, we obtain
\begin{equation}\label{quotient''}
\frac{-\int_{\{\kappa_{\hat\omega}<0\}}\kappa_{\hat\omega}e^{\frac{1}{n}t\psi}\left(\frac{\hat\omega^n}{\omega_0^n}\right)^{\frac1n}\omega(t)^n}{\int_X\omega(t)^n}\ge c_3
\end{equation}
for a positive number $c_3$ depending only on $\omega_0,\delta_1,\delta_2$.
\subsubsection{A lower bound for $\sup_{\{\kappa_{\hat\omega}<0\}}\Phi(t)$.} The above positive lower bound in \eqref{quotient''} in turn gives a lower bound for $\sup_{\{\kappa_{\hat\omega}<0\}}\Phi(t)$. To see this, we integrate \eqref{schlem} and get
\begin{align}
\int_X\omega(t)^n&\ge\int_X(-\kappa_{\hat\omega}+\frac{t}{n})tr_{\omega(t)}\hat\omega\omega(t)^n\nonumber\\
&\ge-\int_{\{\kappa_{\hat\omega}<0\}}\kappa_{\hat\omega}tr_{\omega(t)}\hat\omega\omega(t)^n\nonumber\\
&\ge-n\int_{\{\kappa_{\hat\omega}<0\}}\kappa_{\hat\omega}\left(\frac{\hat\omega^n}{\omega(t)^n}\right)^{1/n}\omega(t)^n\nonumber\\
&=-n\int_{\{\kappa_{\hat\omega}<0\}}\kappa_{\hat\omega}\left(\frac{\omega_0^n}{\omega(t)^n}\right)^{1/n}\left(\frac{\hat\omega^n}{\omega_0^n}\right)^{\frac1n}\omega(t)^n\nonumber\\
&=-n\int_{\{\kappa_{\hat\omega}<0\}}\kappa_{\hat\omega}e^{-\frac1n\Phi(t)+\frac1nt\psi}\left(\frac{\hat\omega^n}{\omega_0^n}\right)^{\frac1n}\omega(t)^n\nonumber\\
&\ge-ne^{-\frac1n\sup_{\{\kappa_{\hat\omega}<0\}}\Phi(t)}\int_{\{\kappa_{\hat\omega}<0\}}\kappa_{\hat\omega}e^{\frac1nt\psi}\left(\frac{\hat\omega^n}{\omega_0^n}\right)^{\frac1n}\omega(t)^n\nonumber,
\end{align}
which, combining with \eqref{quotient''} concludes
\begin{align}
\sup_{\{\kappa_{\hat\omega}<0\}}\Phi(t)\ge n\log nc_3.
\end{align}

\subsubsection{Completion of the proof of Theorem \ref{thm_gap}}
Now we define
$$c_4:=\inf\left\{\int_Xe^v\omega_0^n|v\in PSH(X,(c_0+b_0)\omega_0),n\log nc_3\le\sup_Xv\le\log(c_0+b_0)^n\right\},$$
which is a positive number depending only on $\omega_0,\delta_1,\delta_2$. By \eqref{cm-1} and $\psi\le0$,

\begin{align}\label{cm_2''}
\int_X\left(n\mu_{\hat\omega}\omega_0-Ric(\omega_0)\right)^n&=\lim_{t\to n\mu_{\hat\omega}}\int_X\left(t\omega_0-Ric(\omega_{0})+\sqrt{-1}\partial\bar\partial\Phi(t)\right)^n\nonumber\\
&=\lim_{t\to n\mu_{\hat\omega}}\int_Xe^{\Phi(t)}e^{-t\psi}\omega_{0}^n\nonumber\\
&\ge\liminf_{t\to n\mu_{\hat\omega}}\int_Xe^{\Phi(t)}\omega_{0}^n\nonumber\\
&\ge c_4,
\end{align}
On the other hand, 
\begin{align}\label{cm_2'''}
\int_X\left(n\mu_{\hat\omega}\omega_0-Ric(\omega_0)\right)^n&=\lim_{t\to n\mu_{\hat\omega}}\int_X\left(t\omega_0-Ric(\omega_{0})+\sqrt{-1}\partial\bar\partial\Phi(t)\right)^n\nonumber\\
&\le\int_X\left(2n\hat\epsilon\omega_0-Ric(\omega_0)\right)^n\nonumber\\
&\le\int_X\left(-Ric(\omega_0)\right)^n+c_5\hat\epsilon\nonumber\\
&=\int_X\left(2\pi c_1(K_X)\right)^n+c_5\hat\epsilon,
\end{align}
where $c_5$ is a positive number depending only on the integrals $\int_X\omega_0^s\wedge(-Ric(\omega_0))^{n-s}$, $s=1,...,n$.
Plugging \eqref{cm_2'''} into  \eqref{cm_2''} concludes
\begin{align}
\int_X(2\pi c_1(K_X))^n\ge c_4-c_5\hat\epsilon\nonumber.
\end{align}
Therefore, when $\hat\epsilon\le\min\{\frac{c_4}{2c_5},\frac{c_0}{2n}\}$, we have
\begin{align}
\int_X(2\pi c_1(K_X))^n\ge\frac{c_4}{2}>0\nonumber.
\end{align}

Theorem \ref{thm_gap} is proved.

\subsection{Proof of Theorem \ref{thm_main}}
Note that in the setting of Theorem \ref{thm_main}, the K\"ahler class $[\omega_0]$ has \emph{almost nonpositive holomorphic sectional curvature} (see Definition \ref{defn1}(3)). Then by first-named author's previous work $K_X$ is nef (\cite[Theorem 1.6]{Zys1}) and $X$ does not contain any rational curves (\cite[Theorem 1.10]{Zys1} or \cite[Corollary 4.1, Remark 4.1(5)]{Zys2}). Moreover, by the assumptions of Theorem \ref{thm_main} we may assume, up to passing to a subsequence, that
\begin{equation}\label{ass'n3} 
\lim_{i\to\infty}\mathscr H(\hat\omega_i,\lambda_0;\omega_0)=:\delta_0>0
\end{equation}
Then applying Theorem \ref{thm_gap} gives that

\begin{equation}\label{big_k_x}
\int_Xc_1(K_X)^n>0.
\end{equation}
Therefore, given nefness of $K_X$ and \eqref{big_k_x}, $K_X$ is big and hence $X$ is projective, as pointed out in \cite[Section 2]{DT}. Finally, since $X$ contains no rational curves, by \cite{Ka} (see also \cite[Lemma 5]{WY1} or \cite[Lemma 2.1]{DT}), $K_X$ is ample.

Theorem \ref{thm_main} is proved.

\subsection{Proof of Theorem \ref{thm_main'}} The proof is similar to that of Theorem \ref{thm_gap}. Let's descripe some details. Firstly let $b_0,c_0$ be the same as above (see \eqref{riclower} and \eqref{phiupper}), and define
\begin{align}
c_1:=\sup\{-(b_0+c_0)^n\int_Xu\omega_0^n|u\in PSH(X,(b_0+c_0)\omega_0),\sup_Xu=0\}<+\infty.
\end{align}
which depends only on $\omega_0$. We claim that the number $K_0:=e^{nc_1}$ will satisfy the required property, which can be checked as follows. 

By the assumption we may fix a $c_6>K_0$ such that for every sufficiently large $i$,
\begin{equation}\label{ineq-c6}
\int_{V(\hat\omega_i,\lambda;\omega_0)}\left(\log\frac{(-\kappa_{\hat\omega_i}\hat\omega_i)^n}{\omega_0^n}\right)\omega_0^n\ge c_6.
\end{equation}

For each $i$, we solve the continuity method
\begin{equation}\label{cm2}
(t\hat\omega_i-Ric(\omega_0)+\sqrt{-1}\partial\bar\partial\varphi_i(t))^n=e^{\varphi_i(t)}\omega_0^n,
\end{equation}
which, writing $\omega_i(t):=t\hat\omega_i-Ric(\omega_0)+\sqrt{-1}\partial\bar\partial\varphi_i(t)$, is equivalent to
\begin{equation}\label{cm1}
Ric(\omega_i(t))=-\omega(t)+t\hat\omega_i.
\end{equation}
Similar to the above, the continuity method \eqref{cm2} can be smoothly solved for $t\in(n\mu_{\hat\omega_i},+\infty)$. 
Set $t_i:=2n\mu_{\hat\omega_i}$, 
$$\Phi_i:=\varphi_i(t_i)+t_i\psi_i,$$
where $\psi_i\in C^\infty(X)$ satisfies $\hat\omega_i=\omega_0+\sqrt{-1}\partial\bar\partial\psi_i$ and $\sup_X\psi_i=0$, and 
$$\omega_i:=\omega_i(t_i).$$
As in subsection \ref{proof_gapthm}, we know that, for every sufficiently large $i$, 
\begin{align}
\sup_X\Phi_i\le n\log(b_0+c_0),
\end{align}
and
\begin{align}
\sup_{\{\kappa_{\omega_i}<0\}}\Phi_i&\ge n\log\left(n\frac{\left(\int_{V_i}\omega_0^n\right)\exp\left(\frac{1}{n\int_{V_i}\omega_0^n}\int_{V_i}\log\left(\frac{(-\kappa_{\hat\omega_i}\hat\omega_i^n)^n}{\omega_0^n}\right)\omega_0^n-\frac{c_1}{\int_{V_i}\omega_0^n}\right)}{c_0^{-n}C_{c_0\omega_0}}\right)\nonumber\\
&\ge n\log\left(n\frac{\left(\int_{V_i}\omega_0^n\right)\exp\left(\frac{\frac1n\log c_6-c_1}{\int_{V_i}\omega_0^n}\right)}{c_0^{-n}C_{c_0\omega_0}}\right)\nonumber,
\end{align}
where $V_i:=V(\omega_i,\lambda;\omega_0)$.
Consequently,
$$n\log\left(n\frac{\left(\int_{V_i}\omega_0^n\right)\exp\left(\frac{\frac1n\log c_6-c_1}{\int_{V_i}\omega_0^n}\right)}{c_0^{-n}C_{c_0\omega_0}}\right)\le n\log(b_0+c_0),$$
from which, since $\frac1n\log c_6-c_1>0$, we immediately conclude that 
\begin{equation}\label{volbd-i}
\liminf_{i\to\infty}\int_{V_i}\omega_0^n>0.
\end{equation}
Finally, \eqref{volbd-i} implies $\liminf_{i\to\infty}\mathscr H(\omega_i,\lambda;\omega_0)>0$, then we can apply Theorem \ref{thm_main} to finish the proof.

Theorem \ref{thm_main'} is proved.

\section*{Acknowledgements}
The first-named author thanks professor Huai-Dong Cao for valuable discussions and suggestions on related topics in fall 2018. Both authors are grateful to professors Huai-Dong Cao, Simone Diverio and Xiaokui Yang for their interest in the results, and in particular to professor Valentino Tosatti for valuable comments and criticisms on a previous version of this paper.

\end{document}